\documentclass{article}

\usepackage{amsmath,amssymb,upgreek}

\newtheorem{theorem}{Theorem}[section]

\newtheorem{definition}{Definition}[section]

\newtheorem{proposition}[definition]{Proposition}
\newtheorem{corollary}[definition]{Corollary}

\DeclareMathAlphabet\mathbit
    \encodingdefault\rmdefault\bfdefault\itdefault
\DeclareOldFontCommand{\bi}{\normalfont\bfseries\itshape}{\mathbit}

\newcommand{\be}{\begin{equation}}
\newcommand{\ee}{\end{equation}}

\def\fakebold#1{\relax\ifvmode\leavevmode\fi%
\ifmmode%
\setbox0=\hbox{$#1$}%
\else%
\setbox0=\hbox{#1}%
\fi%
\kern-.02em\copy0 \kern-\wd0%
\kern .04em\copy0 \kern-\wd0%
\kern-.0125em\raise.02em\box0%
}%

\renewcommand{\geq}{\geqslant}
\renewcommand{\leq}{\leqslant}
\newcommand{\dd}{\mathrm{d}}
\newcommand{\e}{\mathrm{e}}

\newcommand{\A}{\mathbf{A}}
\newcommand{\B}{\mathbf{B}}
\newcommand{\C}{\mathbf{C}}
\newcommand{\G}{\mathbf{G}}
\renewcommand{\H}{\mathbf{H}}
\newcommand{\K}{\mathbf{K}}

\newcommand{\N}{\mathbf{N}}
\newcommand{\R}{\mathbf{R}}

\newcommand{\X}{\mathbf{X}}
\newcommand{\I}{\mathbf{I}}
\newcommand{\F}{\mathbf{F}}
\newcommand{\V}{\mathbf{v}}
\newcommand{\w}{\mathbf{w}}
\newcommand{\D}{\mathrm{D}}
\newcommand{\Epsilon}{E}

\begin{document}

\bibliographystyle{plain}

\title{On solutions of matrix-valued convolution equations, anisotropic fractional derivatives and their applications in linear and non-linear anisotropic viscoelasticity\footnote{An new version of [1] with improved notation and Theoorem 3.2.}}
\author{Andrzej Hanyga\\
ul. Bitwy Warszawskiej 14/52\\
\rm{02-366} Warszawa, Poland}

\date{}

\maketitle 

\vspace{0.3cm}

%\hfil{Dedicated to the memory of Ma{\l}gorzata Anna Seredy\'{n}ska\\ \vspace{0.5cm}}

\noindent\textbf{Notation.}\\
\begin{enumerate}
\item $\D = \dd/\dd t$
\item $\mathbb{R}$ - the set of real numbers
\item $\mathbb{N}$ - the set of non-negative integers 
\item $[a,b[$ - the set of $x$ such that $a \leq x < b$
\item $\mathcal{S}$ - the set of symmetric matrices of rank N
\item $\mathcal{S}_+$ - the set of positive semi-definite elements of $\mathcal{S}$  
\end{enumerate}
\vspace{0.3cm}

\begin{abstract}
A relation between matrix-valued complete Bernstein functions and matrix-valued Stieltjes functions is applied to prove that the
solutions of 
matricial convolution equations 
with extended LICM kernels belong to special classes of functions. In particular the cases of the solutions of the 
viscoelastic duality relation 
and the solutions of the matricial Sonine equation are discussed, with applications in anisotropic linear viscoelasticity 
and a generalization of fractional calculus. 

In the first case it is in particular shown that duality of completely monotone relaxation functions and Bernstein creep functions
in general requires inclusion in the relaxation function of a Newtonian viscosity term in addition to the memory effects represented 
by the completely monotone kernel.

We define anisotropic generalized fractional derivatives (GFD)
by replacing the kernel $t^{-\alpha}/\Gamma(1-\alpha)$ of the Caputo derivatives with completely monotone matrix-valued kernels 
which are weakly singular at 0.
\end{abstract}


\section{Introduction}

A recent idea \cite{HanDualityNew} of simplifying the proof of a viscoelastic duality relation between the 
tensor-valued relaxation modulus and the tensor-valued creep function \cite{HanDuality}, based on a
relation between matrix-valued complete Bernstein functions and matrix-valued Stieltjes functions 
\cite{HanAnisoWaves} has led me to consider a more 
general application of this method to convolutions of symmetric matrix-valued functions. 

In this reference I have decided to base the analysis of the viscoelastic duality relation on the analysis of the solutions of 
the matrix-valued solutions of the Sonine equation \cite{SamkoKilbasMarichev}.

In the case of the Sonine equation it is assumed here that one of the functions is locally integrable near zero and completely monotone (LICM). 
The Sonine equation was examined in much detail in \cite{SamkoCardoso}, but the authors did not assume that one of the functions
was LICM. They constructed the inverse operator for the convolution equation $k(t)\ast x(t) = f(t)$. The inverse operator however
involves the solution $l(t)$ of the Sonine equation $k(t) \ast l(t) = 1$. Existence of such a function is the subject of our investigation 
and we prove it for the kernel $k$ in the extended LICM class, defined below. The last-mentioned problem is also studied in 
\cite{KochubeiArxiv} for real-valued functions but we consider more general symmetric matrix-valued functions.

I reexamine Kochubei's suggestion \cite{KochubeiArxiv} that the solutions of the Sonine equation could be used to construct a 
generalization of the concepts of a derivative and integration operators along the lines of fractional calculus. I however allow for 
matrix-valued LICM kernels which introduce 3D modeling of anisotropic effects in the context of non-linear relaxation
equations. 

I also show that in order to satisfy a Sonine equation the LICM function $k(t)$ must be singular. If $k(t)$ is a singular
LICM function, then it satisfies the Sonine equation with an associated function $l(t)$, which is also a singular LICM function.

The results presented here are relevant for anisotropic linear and non-linear viscoelasticity. They also open new methods 
for dealing with problems involving memory effects.

\section{Convolutions equations for matrix-valued functions.} \label{First}

The convolution of two square matrix-valued functions $\F$ and $\G$ of the same rank $N$ is defined by the equation
\begin{equation}
(\F \ast \G)(t) := \int_0^t \F(t-s) \, \G(s)\, \dd s 
\end{equation}

The Laplace transform $\tilde{\F}(p)$ of a matrix-valued function $\F(t)$ is defined as usual by the formula
$$\tilde{\F}(p) := \int_0^\infty \e^{-p t} \, \F(t)\, \dd t$$
for every $p \in \mathbb{C}$ such that the integral exists. 

It is easy to check the identity
\begin{equation} \label{LaplId}
\int_0^\infty \e^{-p t}\, \F(t)\ast\G(t) \, \dd t = \tilde{\F}(p)\, \tilde{\G}(p), \; p > 0
\end{equation}
for arbitrary matrix-valued functions defined on $[0,\infty[$ provided both Laplace transforms on the right-hand side exist.

We shall equip the convolution algebra with a unit element $\mathrm{U}$:
\begin{equation} \label{unity}
\mathrm{U}\, f = f
\end{equation}
(see Appendix~A for details).

The unit operator $U$ is a convolution with a Borel measure $\upsilon$ on $[0,\infty[$ defined by the formula 
$\int_{[0,\infty]}  f(r)\, \upsilon(\dd r) = f(0)$ for every continuous function $f$ on $[0,\infty[$. Indeed, 
\begin{equation} \label{unit}
(U f)(r) = \int_{[0,\infty[} f(r - s)\, \, \upsilon(\dd s) = f(r).
\end{equation}
The unit operator $U$ is a convolution with a Borel measure $\upsilon$ on $[0,\infty[$ defined by the formula 
$\int_{[0,\infty]}  f(r)\, \upsilon(\dd r) = f(0)$ for every continuous function $f$ on $[0,\infty[$. Indeed, 
$$ (U f)(r) = \int_{[0,\infty[} f(r - s)\, \, \upsilon(\dd s) = f(r).$$

Extending identity~\eqref{LaplId} to \eqref{unity} we have
$\tilde{\upsilon}(p) \tilde{f}(p) = \tilde{f}(p)$, whence
\begin{equation}
\tilde{\upsilon}(p) = 1, \; p \geq 0
\end{equation}

Consider the general convolutional equation 
\begin{equation} \label{def}
\A_1\, \X(t) + \F(t)\ast\X(t) = \R(t)
\end{equation}
where $\F(t)$ and $\R(t)$ are two square matrix-valued functions defined for $t \in\, ]0,\infty[$ in a class to be specified, 
$\A_1$ is a positive semi-definite symmetric matrix (possibly zero), while 
$\X(t)$ is a square matrix-valued function defined by equation~\eqref{def}. We shall examine the properties of the function $\X(t)$.
The ranks of the matrices are equal and will be denoted by $N$. Let $\mathcal{S}$ denote the space of symmetric real square matrices 
of a fixed rank $N$. 

If $\A_1 = 0$ and $\R(t)$ is unit matrix for $t \geq 0$, then equation~\eqref{def} is a generalization of the {\em Sonine equation} 
\cite{SamkoCardoso} for matrix-valued functions.

\begin{definition}
A matrix-valued function $\F:\, ]0,\infty[ \rightarrow \mathcal{S}$ is said to be {\em completely monotone} (CM) if it is 
infinitely differentiable and for every vector
$\V \in \mathbb{R}^N$ the following inequalities are satisfied
\begin{equation} \label{CM}
\forall t > 0 \; \forall n \in \mathbb{N} \; \; (-1)^n \, \D^n \V^\mathsf{T}\, \F(t)\, \V \geq 0 
\end{equation}
\end{definition}
The above definition allows for a singularity at 0. 

\begin{definition}
A matrix-valued function $\F:\; ]0,\infty[ \rightarrow \mathcal{S}$ is said to be {\em locally integrable completely monotone} 
(LICM) if it is CM and integrable over $]0,1]$.
\end{definition}

\begin{theorem} \cite{HanDuality}\\
Every matrix-valued LICM function $\F$ can be expressed in the form
\begin{equation} \label{LICM}
\F(t) = \int_{[0,\infty[} \e^{-r t} \, \H(r) \, \mu(\dd r),  \; t > 0
\end{equation}
where $\mu$ is a Borel measure on $[0,\infty[$ satisfying the inequality
\begin{equation} \label{ineq}
\int_{[0,\infty[} (1 + r)^{-1}\, \mu(\dd r) < \infty 
\end{equation}
and $\H(r)$ is a matrix-valued function on $[0,\infty[$ satisfying the bound 
$\vert \H(r) \vert \leq 1$ except perhaps on a subset of $[0,\infty[$ of $\mu$-measure zero. 
\end{theorem}

\begin{proposition} \label{decrease}
If the matrix-valued function $\F$ is LICM, then the function $\tilde{\F}(p)$ tends to 0 for $p \rightarrow \infty$.
\end{proposition}
Note that this result also implies that every LICM function $\F$ satisfies the identity
$$\lim_{t\rightarrow 0} \int_0^t \F(s) \, \dd s = 0.$$

\noindent\textbf{Proof.}\\
For $p > 1$ 
$$\vert \tilde{\F}(p) \vert = \int_{[0,\infty[} (p + r)^{-1}\, \mu(\dd r) < \int_{[0,\infty[} (1 + r)^{-1}\, \mu(\dd r) < \infty,$$
hence the thesis follows from the Lebesgue Dominated Convergence Theorem.\\

\mbox{ }\hfill$\Box$

\begin{definition}
The set $\mathcal{E}$ consists of matrix-valued Borel measures of the form
$$\Phi = \upsilon \, \A + \F\, \dd t$$
where $\A$ is a positive semi-definite symmetric matrix, $\F$ is an LICM $\mathcal{S}$-valued function and $\dd t$ is the Lebesgue measure
on $[0,\infty[$. 
\end{definition}

The set $\mathcal{E}$ appears in a natural way in the context of solutions of convolution equations. It also allows for a Newtonian viscosity component in viscoelastic relaxation functions (Section~\ref{linear}).

\begin{definition}
A function $\B:\; ]0,\infty[ \rightarrow \mathcal{S}$ is said to be a {\em Bernstein function} (BF) if $\B(t)$
is positive semi-definite and differentiable for $t > 0$ and its derivative is CM. 

We denote by $\mathcal{B}$ is the set of all the Bernstein functions.
\end{definition}

If $\B$ is an $\mathcal{S}$-valued Bernstein function, then for every $\V \in \mathbb{R}^N$ 
the function $t \rightarrow \V^\mathsf{T}\, \B(t)\, 
\V$, $\;t > 0$, is non-decreasing and continuous, hence it has a finite limit at $t = 0$. Hence the limit 
$\lim_{t\rightarrow 0} \B(t)$ exists. We shall therefore consider Bernstein functions as defined on $[0,\infty[$ with 
$\B(0) = \lim_{t\rightarrow 0} \B(t)$. 

It follows that the derivative $\B^\prime$ of $\B$ on $]0,\infty[$ is a completely monotone function. It is also locally 
integrable because 
$$\int_0^1 \B^\prime(t)\, \dd t = \B(1) - \B(0), \; \B^\prime(t) \geq 0,$$ 
hence it is LICM.

A general Bernstein 
function is obtained by integrating the measure $\upsilon\, \B_0 + \F \, \dd t$ over intervals $[0,t]$ for a symmetric and 
positive semi-definite matrix $\B_0$ and an LICM matrix-valued function $\F$ .

An $\mathcal{S}$-valued LICM function $\F$ is locally integrable and non-increasing, hence its Laplace transform 
$\tilde{\F}(p)$ is defined for all $p > 0$. 

An $\mathcal{S}$-valued Bernstein function $\B$ also has a Laplace transform
defined for $p > 0$ because $$\int_0^\infty \e^{-p t}  \B(t) \, \dd t =  
\int_0^\infty \e^{-p t}\, \int_0^t \F(s)\, \dd s \,\dd t = \frac{1}{p} \tilde{\F}(p),$$
where $\F$ is the LICM derivative of $\B$.

Let $\mathcal{S}_+$ denote the set of positive semi-definite symmetric $N \times N$ matrices.

\begin{proposition} \label{PCM}
If $\F$ is a symmetric matrix-valued LICM function and for each vector $\V \in \mathbb{R}^N$ the function 
$\V^\mathsf{T} \, \F(t)\, \V, t > 0$, is not identically zero, then 
the matrix $\tilde{\F}(p)$ is invertible for every $p > 0$.
\end{proposition}

\noindent\textbf{Proof.}\\
For each non-zero $\V \in \mathbb{R}^N$ there is a real $t_1(\V) > 0$ such that $\V^\mathsf{T}\, \F(t_1(\V)) \, \V > 0$. The function
$\V^\mathsf{T}\, \F(t)\, \V$ is continuous, hence it is positive on some interval $\mathcal{I} \subset ]0,\infty[$, while it is 
non-negative on $]0,\infty[$. Hence for every non-zero $\V \in \mathbb{R}^N$, $p > 0$ we have $\V^\mathsf{T} \,\tilde{\A}(p)\, \V > 0$. 
The matrix $\tilde{\F}(p)$ is symmetric and positive definite, hence it is invertible.\\
\mbox{ }\hfill$\Box$

\begin{proposition} \label{PBF}
If $\B$ is a symmetric matrix-valued Bernstein function and for each vector $\V \in \mathbb{R}^N$ the function 
$\V^\mathsf{T} \, \B(t)\, \V, t > 0$, is not identically zero, then 
the matrix $\tilde{\B}(p)$ is invertible for every $p > 0$.
\end{proposition}

The proof of Proposition~\ref{PBF} is analogous to Proposition~\ref{PCM}.

Equation~\eqref{LaplId} is equivalent to the equation
\begin{equation} \label{Lap}
\left[\A_1 + \tilde{\F}(p) \right]\, \tilde{\X}(p) = \tilde{\R}(p)
\end{equation}

If either the matrix $\A_1$ is positive definite or the matrix $\tilde{\F}(p)$ is invertible for $p > 0$, 
then the matrix $\A_1 + \tilde{\F}(p)$ is invertible for $p > 0$. In this case 
the unique solution of \eqref{Lap} is 
\begin{equation} \label{Lap1}
\tilde{\X}(p) = \left[ \A_1 + \tilde{\F}(p)\right]^{-1}\, \tilde{\R}(p), \;\; p > 0
\end{equation}
This is true in particular if $\F$ is an $\mathcal{S}_+$-valued LICM function. 

For $\R(t) = \I$ ($t > 0$), where $\I$ is the identity operator on $\mathbb{R}^N$, $\tilde{\R}(p) = p^{-1}\, 
\I$ and $\tilde{\X}(p) = \left[ p \, \A_1 + p\, \tilde{\F}(p) \right]^{-1}$.
We shall show that in this case $\X$ is an LICM function. 

For $\R(t) = t \, \I, t > 0$, we have $\tilde{\R}(p) = p^{-2}\, \I$ and
\begin{equation}
p \, \tilde{\X}(p) = \left[ p\, \A_1 + p\, \tilde{\F}(p) \right]^{-1}
\end{equation}
We shall show that under these assumptions the solution $\X$ of equation~\eqref{def} is an $\mathcal{S}_+$-valued Bernstein function.

Similarly, if $\A$ is an $\mathcal{S}_+$-valued BF and $\R(t) = t \, \I, t > 0$, then $\X$ is an $\mathcal{S}_+$-valued LICM function.

By transposition the results obtained below also apply to equations of the form
$$\X(t)\ast (\A_1 + \F(t)) = \R(t).$$

\section{Main theorem.}

Assume that $\A_1$ is positive semi-definite symmetric matrix of rank $N$.

\begin{theorem} \label{thm3}
If either\\ (1) $\A_1$ is a positive definite symmetric matrix\\ or else\\ (2)
$\F$ is a non-zero $\mathcal{S}_+$-valued LICM function, satisfying the conditions \\
\parbox{\textwidth}{$(\ast)$ for each non-zero vector $\V \in \mathbb{R}^N$ the function 
$\V^\mathsf{T} \, \F(t)\, \V, t > 0$, is not identically zero,}   
and \\
\parbox{\textwidth}{$(\ast\ast)$ the limit 
\begin{equation} \label{zz}
\F_0 := \lim_{t\rightarrow 0} \F(t)^{-1}
\end{equation}
 exists,} \\
then equation~\eqref{def} with $\R(t) = \I$
has a unique solution $\X(t)$, where $\X$ is an LICM function in the first case and 
$\X(t) = \upsilon(t)\, \F_0  + \G(t)$,  where $\G$ is an $\mathcal{S}_+$-valued LICM function, in the second case. 
\end{theorem}

In general the solution of equation~\eqref{def} $\X \in \mathcal{E}$.\\

\noindent\textbf{Remark.}
Concerning Condition~$(\ast\ast)$, we begin with the remark that the matrix $\F(t)$ is invertible for sufficiently large $t$. 
Indeed, on account of 
Condition~$(\ast)$ each of its eigenvalues $a_n(t)$ 
($n = 1, \ldots N$) is a CM function and is positive for some $t_n$. Consequently the matrix $\F(t)$ is invertible 
for $t > \sup\{t_n \mid n = 1, \ldots N\}$. Hence it makes sense to inquire whether the limit $\lim_{t\rightarrow 0} \F(t)^{-1}$ exists.\\

\noindent\textbf{Remark.}\\
The equation $k\ast l = 1$ in for locally integrable real-valued functions $k$ and $l$ is known as the Sonine equation 
\cite{SamkoKilbasMarichev}. If for a given $k \in \mathcal{L}^1_{\mathrm{loc}}([0,\infty[)$ there is an 
$l \in \mathcal{L}^1_{\mathrm{loc}}([0,\infty[)$ satisfying the above equation, then $k$ is called a {\em Sonine function}, 
while $k, l$ are known as a {\em Sonine pair}. Sonine pairs are studied in some detail in \cite{SamkoCardoso}. 
Theorem~\ref{thm3} asserts in particular that every LICM function or matrix-valued LICM function is a Sonine function and 
in this case the Sonine pair consists of two LICM functions.
For real-valued functions this fact has apparently been discovered by Kochubei \cite{KochubeiArxiv}.

However, not every Sonine pair consists of LICM functions. A counterexample is the Sonine pair $k_\lambda(t) := t^{-\lambda/2}\,
J_{-\lambda}(2 t^{1/2})$ with the Laplace transform $\tilde{k}_\lambda(p) = \exp(-1/p) \, p^{\lambda-1}$ 
(\cite{BatemanProject} p. 185 (30)) and
$l_\lambda(t) t^{(\lambda-1)/2} \, I_{\lambda-1}(2 t^{1/2})$ with $\tilde{l}_\lambda(p) = \exp(1/p)\, p^{-\lambda-2}$
(\cite{BatemanProject} p. 197 (18))
for $\lambda > 0$. The function $k_\lambda$ changes sign and therefore is not CM, for example
$k_{1/2}(t) = \sqrt{2/\uppi} \, \cos(2 t^{1/2})/t^{3/4}$. 

The following matrix-valued Sonine pairs are of particular interest:
\begin{enumerate}
\item $k(t) \, \K_0$ and $l(t)\, \K_0^{\;-1}$, where $k, l$ are a Sonine pair of LICM functions;
\item $\mathrm{diag} \{k_n(t), n=1,... m\}$ and $\mathrm{diag} \{l_n(t), n=1,... m\}$, where 
$(k_n, l_n)$ are Sonine pairs of LICM functions for $n = 1, \ldots, m$.
\end{enumerate}

Many LICM functions are known \cite{MillerSamko01}, but it is often more difficult to find the other member of the Sonine pair.
The simplest Sonine pair of CM functions is $k(t) = t^{\alpha-1}/\Gamma(\alpha)$ and $l(t) = t^{-\alpha}/\Gamma(1-\alpha)$,
$0 < \alpha < 1$. Using the Laplace transforms $\tilde{k}(p) = (p + \lambda)^{-\alpha}$ and  
$\mathcal{L}[\Gamma(-\alpha,\lambda t)](p) = \Gamma(-\alpha)\, \lambda^{-\alpha} \left[\lambda^\alpha - (\lambda + p)^\alpha\right]/p$
one gets another pair  $k(t) = t^{\alpha-1}\, \e^{-\lambda t}/\Gamma(\alpha)$, 
$\lambda > 0$, with $l(t) = \lambda^\alpha\,\left[1 - \Gamma(-\alpha,\lambda t)/\Gamma(-\alpha)\right]$,
$\lambda \geq 0$, $0 < \alpha < 1$. 

It is also interesting that for an arbitrary analytic function $k(t)$ there is another analytic function $l(t)$ such that
$k(t)\, t^{\alpha-1}/\Gamma(\alpha)$ and $l(t) \, t^{-\alpha}/\Gamma(1-\alpha)$ are a Sonine pair and there is an algorithm for 
calculating the power series of $l(t)$ given the power series for $k(t)$ \cite{SamkoCardoso, Wick}.\\
\vspace{0.2cm}

\noindent\textbf{Proof of Theorem~\ref{thm3}.}

The Laplace transform $\tilde{\F}(p)$ is a symmetric positive definite matrix for every  $p > 0$. For $p > 0$ 
equation~\eqref{def} is equivalent   equation $\tilde{\mathbf{X}}(p) = \left[p\, (\A_1 + \tilde{\F}(p))\right]^{-1}$. The inverse on 
the right-hand side
exists for $p > 0$ if $\A_1$ is positive definite or else in view of Condition~($\ast$) and Proposition~\ref{PCM}. 
The right-hand side of the last equation is the algebraic inverse of a matrix-valued CBF, hence it is a matrix-valued Stieltjes 
function of the form 
\begin{equation}  \label{eq5}
\mathbf{K}(p) = \mathbf{C} + \int_{[0,\infty[} (p + r)^{-1}\, \mathbf{H}(r)\, \mu(\dd r).
\end{equation}
where $\mathbf{C} \geq 0$, $\mathbf{H}$ is a measurable symmetric matrix-valued function bounded $\mu$-almost everywhere 
and $\mu$ is a Borel measure satisfying inequality~\eqref{ineq} (Theorem~\ref{inv}). The second term on the right-hand side of equation~\eqref{eq5} 
is a Laplace transform $\tilde{\G}(p)$ of the LICM matrix-valued function 
\begin{equation} \label{FLICM}
\G(t) := \int_{[0,\infty[} \e^{-r t}\, \mathbf{H}(r) \, \mu(\dd r)
\end{equation}

By Proposition~\ref{decrease} 
$$\mathbf{C} = \lim_{p \rightarrow \infty} \mathbf{K}(p) = 
\left[ \lim_{p \rightarrow \infty}\,(p\, \A_1 + p\, \tilde{\F}(p))\right]^{-1} $$ 

Again by Proposition~\ref{decrease}, 
if $\A_1 > 0$, then $\mathbf{C} = 0$. If $\A_1 = 0$, then in view of Condition~$(\ast\ast)$ $\mathbf{C} = \F_0$.

Applying the inverse Laplace transformation to \eqref{eq5} we conclude that $\X = \G$ if $\A_1 > 0$, $\X(t) = 
\upsilon(t)\, \F_0 + \G(t)$ if $\A_1 = 0$.

Uniqueness of the solution $\X(p)$ follows from the fact that in view of the invertibility of the matrix $\A_1 + \tilde{\F}(p)$ 
the equation $(\A_1 + \tilde{\F}(p)) \, \tilde{\X}(p) = 0$ for $p > 0$ implies that $\tilde{\X}(p) = 0$ for $p  \geq 0$.\\

\mbox{ }\hfill$\Box$

\begin{theorem}\label{thm30}
If $\A_1 = 0$ and $\F$ is a singular LICM matrix-valued function, then the solution $\X$ of equation~\eqref{def} is a singular LICM function $\G$.
\end{theorem}

\noindent\textbf{Proof.}\\
Since $\A_1 = 0$, the solution of equation~\eqref{def} is an LICM matrix-valued function $\G$ and
$$\lim_{t\rightarrow 0} \G(t)^{-1} = \lim_{p\rightarrow\infty} \left(p \,\tilde{\G}(p) \right)^{-1} = 
\lim_{p\rightarrow\infty} \tilde{\F}(p) = \lim_{t\rightarrow 0} \int_0^t \F(s)\, \dd s= 0, $$
because the function $\F$ is locally integrable.
\mbox{ }\hfill$\Box$

\begin{corollary} \label{cor4}
For $n \in \mathbb{N}$ the convolution equation \eqref{def} with $\R(t) = (t^n/n!) \, \I$  and an LICM matrix-valued 
function $\F$ satisfying Condition~$(\ast)$ has a unique solution $\X$ which is an $n$-fold indefinite integral of an 
element of $\mathcal{E}$.

In particular, for $n = 1$ the solution $\X$ is a matrix-valued Bernstein function and $\X(0) = 0$ if $\A_1 > 0$, while 
$\X(0) = \F_0$ if $\A_1 = 0$, with $\F_0$ given by equation~\eqref{zz}. 
\end{corollary}

For the proof of the last statement note that $\lim_{t \rightarrow 0} \int_0^t \F(s)\, \dd s = 0$ because $\F$ is integrable
in a neighborhood of 0.

Here is a complementary result for the last statement:
\begin{corollary} \label{cor5}
Let $\B$ be an $\mathcal{S}_+$-valued Bernstein function satisfying the condition\\ \vspace{0.2cm}
\parbox{\textwidth}{$(\ast\ast\ast)$ For every $\V \in \mathbb{R}^N$ the function $\V^\mathsf{T}\, \B(t)\, \V$, $t > 0$,  is not constant.}\\
The equation
\begin{equation} \label{Beq}
\B \ast \X = t \, \I
\end{equation}
has a unique solution $\X$ and this solution is an $\mathcal{S}$-valued LICM function. 
\end{corollary}
\noindent\textbf{Proof.}\\

Differentiating equation~\eqref{Beq} with respect to $t$ one gets equation~\eqref{def} with $\R(t) = \I$ and $\F = \B^\prime$, an LICM function
satisfying Condition~$(\ast)$.
The thesis then follows from Theorem~\ref{thm3}. 

\mbox{ }\hfill$\Box$

\section{Anisotropic generalized fractional derivatives (GFD).}

The term "derivative" has been improperly applied to fractional derivatives although they do not satisfy the Leibniz property, 
which is part of the definition of a derivative. Fractional derivatives have however provided useful tools for constructing equations
which in some
sense interpolate between differential equations of varying orders. We shall now show that similar operators can be constructed
for a kind of anisotropic generalized fractional derivative (GFD) which might be useful to construct anisotropic relaxation 
equations and for other purposes.

The first application of this result is the solution of a convolution equation
\begin{equation}
\F(t)\ast \V(t) = \mathbf{f}(t)
\end{equation}
where $\F$ is a given matrix-valued function.
Since 
\begin{equation} \label{zi}
\X(t)\ast \F(t) = \I, 
\end{equation}
we have
$$\int_0^t \V(t)\, \dd t  = \X(t)\ast \mathbf{f}(t)$$
and therefore
\begin{equation}
\V(t) = \D [\X\ast \mathbf{f}(t)] 
\end{equation}

The LICM matrix-valued function $\F(t)$ can have a singularity at 0 such that for every $\V \in \mathbb{R}^N$ the limit 
$\lim_{t \rightarrow 0} \, \V^\mathsf{T}\,\F(t)\, \V = \infty$. In this case $\F_0 = 0$, where $\F_0$ is defined by 
equation~\eqref{zz}. We shall say that the function $\F$ is {\em singular} if $\F_0 = 0$.

If $\F$ is a singular LICM matrix-valued function, then by Theorem~\ref{thm3} the solution $\X = \G$ of the Sonine equation $\F \ast \X = \I$  is an LICM function and we can define
the $\F$-{\em derivative} by the formula
\begin{equation} \label{diffdef}
\D_\F \, \V(t) = \D \left[ \F(t)\ast \V(t) \right] - \F(t)\, \V(0) \equiv \F \ast \D \V
\end{equation}
for every  absolutely continuous function $\V:\; [0,\infty[ \rightarrow \mathbb{R}^N$.
 If $\X = \G$ (a singular LICM function) is the solution 
of the convolution equation~\eqref{def} with $\R(t) = \I$, then the  $\F$-integral operator is defined by the formula
\begin{equation}
\mathrm{J}_\F \, \V(t) := \G\ast\V(t).
\end{equation}

It follows from Theorem~\ref{thm30} that the function $\G$ is singular. 

We then have 
\begin{theorem} \label{thm4}
Let $\F$ be a singular LICM matrix-valued function.

The following relations hold 
\begin{eqnarray}
\mathrm{J}_\F\, \D_\F \, \w(t) = \w(t) - \w(0) \;\; \mathrm{for}\;\; \w \in AC([0,\infty[)\\
\D_\F\, \mathrm{J}_\F \, \V(t) = \V(t) \; \;\mathrm{for}\;\; \V \in \mathcal{L}^1_\mathrm{loc}([0,\infty[)
\end{eqnarray}
\end{theorem}
\noindent\textbf{Proof}\\
(1) The identity $\G \ast \F = \mathbf{1}$ implies that 
$$(\mathrm{J}_\F \, \D_\F \, \w)(t) = \int_0^t \G(s) \left. \frac{\dd}{\dd \tau} \int_0^\tau \F(\tau - r) \,\w(r)\, \dd r 
\right|_{\tau=t-s} \, \dd s - \w(0)$$
The first term on the right-hand side equals
$$\int_0^t \G(s) \frac{\dd}{\dd t} \int_0^{t-s} \F(t -s - r) \, \w(r)\, \dd r\, \dd s = \frac{\dd}{\dd t} (\G\ast \F \ast \w)(t) =
\frac{\dd}{\dd t} (\mathbf{I}\ast \w)(t) = \w(t)$$
q.e.d.

\noindent (2) Let $\w := \G \ast \V$.
On account of the identity $\F \ast \G = \mathbf{I}$ 
$$\D_\F\, \mathrm{J}_\F\, \V = \frac{\dd}{\dd t} (\F\ast \G\ast \V)(t) - \F(t)\, \w(0) = \V(t) - \F(t)\, \w(0)$$

It remains to prove that $\w(0) = 0$. $\G$ is an LICM function, hence it has the form \eqref{FLICM} with $\mu$ satisfying 
equation~\eqref{ineq} and 
$\vert \H(r)\, \vert \leq 1$. Hence 
$$\vert \w(t) \vert^2 \leq \left[ \int_0^t \int_{[0,\infty[} \e^{-r s} \, \mu(\dd r)\, \dd s\right] \,\int_0^t |\V(t - s)| \, \dd s$$

For $t \leq 1$ the second factor is bounded from above by a constant 
$$\int_0^1 | \V(s)| \, \dd s  < \infty$$
because $\V$ is assumed locally integrable.
The first factor equals
\begin{equation} \label{f1}
\int_{[0,\infty[} \frac{1 - \e^{-r t}}{r} \, \mu(\dd r)
\end{equation}
From the inequality $\e^x - 1 \leq x\, \e^x$ ($x \geq 0$) follows  the inequality $1 - \e^{-x} \leq x$. We shall apply this 
inequality for $r \in [0,\infty[$, noting that 
$\mu([0,1]) < \infty$ because of \eqref{ineq} with the inequality $1 \leq 2/(1 + r)$ valid for $r \leq 1$. For $r > 1$ we shall note that 
$1/r \leq 2/(1 + r)$. Hence expression~\eqref{f1} is bounded by 
$$t  \, \mu([0,1]) + 2 \int_{]1,\infty[} \left( 1 - \e^{-r t}\right)\, (1 + r)^{-1} \, \mu(\dd r),$$
which tends to 0 as $t \rightarrow 0$ on account of \eqref{ineq} and the Lebesgue Dominated Convergence Theorem. Thus $\w(0) = 0$ and the theorem has been proved.

\mbox{ }\hfill$\Box$\\

The new derivative concept provides a new approach to modeling stress relaxation in anisotropic and non-linear viscoelastic media.
A possible relaxation equation could have the form 
\begin{equation} \label{eqrela}
\D_\F \, \Sigma = \mathbf{K}(\Sigma,\Epsilon)
\end{equation}
where $\mathbf{K}$ is a rank-2 tensor-valued function of two rank-2 tensor-valued arguments.

In view of equation~\eqref{diffdef}$_2$ and \eqref{zi} equation~\eqref{eqrela} is equivalent to the equation
\begin{equation}
\Sigma = \Sigma(0) + \X\ast\, \mathbf{K}(\Sigma,\Epsilon)
\end{equation}

In the scalar case Theorem~\ref{thm4} applies only to (weakly) singular CM kernels $\F$ such as $t^{-\alpha}/\Gamma(1-\alpha)$. 

\section{Application to anisotropic linear viscoelasti\-city.}\label{linear}

The relaxation modulus $F_{jklm}(t)$ and the creep function $C_{jklm}(t)$ in 3-dimensional linear viscoelasticity are defined by
 the two constitutive
laws which are assumed equivalent to each other
\begin{eqnarray}
\Sigma_{jk} = N_{jklm} \, \D{\Epsilon}_{lm} + F_{jklm} \,\ast \, \D{\Epsilon}_{lm}, \; j, k = 1,\ldots 3\\
\Epsilon_{jk} = C_{jklm}\, \ast \, \D{\Sigma}_{lm} \; j, k = 1,\ldots 3\
\end{eqnarray}
(summation over repeated indices is assumed)
where the symmetric tensors $\Sigma_{kl}$ and $\Epsilon_{kl}$ denote the stress and strain tensors, respectively.
We assume the usual symmetries
\begin{eqnarray} \label{symmetries}
F_{ijkl} = F_{jikl} = F_{klij} \nonumber \\
C_{ijkl} = C_{jikl} = C_{klij} \\
N_{ijkl} = N_{jikl} = N_{klij} \nonumber \\
\mathrm{for }\;\; i, j , k ,l = 1\, \ldots 3 \nonumber
\end{eqnarray}
Equivalence of the two constitutive equations is ensured by the relation
\begin{equation} \label{convVE}
N_{jklm} \, C_{jklm} + F_{jklm} \ast C_{lmrs} = t \, (\delta_{jr}\, \delta_{ks} + \delta_{js}\, \delta_{kr}),  \;\; j, k, r, s = 1,\ldots, 3,\; t > 0
\end{equation}
The first term on the left-hand side represents a Newton viscosity component. It is assumed to satisfy  the inequalities 
$N_{jklm} \, e_{jk}\, e_{lm} \geq 0$ for every symmetric tensor $e_{kl}$. 
   
Defining the index $I$, $1 \leq I \leq 6$,
$I = k$ for the pair $kk$, $1 \leq k \leq d$ and $I = m$ for the pair $kl$, $k \neq l, m$, $l \neq m$,
 $R_{IJ} = f(I) \, f(J)  \, R_{IIJJ}$, with 
$f(I) = 1$ if $1 \leq I \leq 3$ and $f(I) = \sqrt{2}$ if $4 \leq I \leq 6$. In this notation equation~\eqref{convVE} assumes the 
form of a matrix convolution of two symmetric 6 $\times$ 6 matrix-valued functions 
\begin{equation} \label{anisorel-t}
N_{IJ} \, C_{JK} + F_{IJ}\, \ast C_{JK} = t \, \delta_{JK}  \;\; J, K =1,\ldots 6
\end{equation}
In the Laplace domain equation~\eqref{anisorel-t} is equivalent to
\begin{equation} \label{anisorel1}
N_{IJ} \, \widetilde{C_{JK}} + \widetilde{F_{IJ}}\, \widetilde{C_{JK}} = p^{-2} \, \delta_{JK}, \;\; J, K =1,\ldots 6
\end{equation}

The matrix-valued function $F_{IJ}$ is LICM if for every $\mathbf{v} \in \mathbb{R}^6$ the function 
$t \rightarrow v_I \, v_J \, F_{IJ}(t)$ is LICM.
The above statement is equivalent to $e_{ij} \, e_{kl} \, F_{ijkl}(t)$ being LICM for every symmetric 3 $\times$ 3 matrix $e_{ij}$. 

We shall use the notation
$\N$, $\mathbf{R}$, $\mathbf{C}$ for the 6 $\times$ 6 matrices $\N_{IJ}$,  $R_{IJ}$ and $C_{IJ}$. The inequalities  
$\N \geq 0$, $\N > 0$ are short-hand for the inequalities $N_{ijkl}\, e_{ij}\, e_{kl} \geq 0$ for an arbitrary rank-2 tensor $e_{kl}$ and
$N_{ijkl}\, e_{ij}\, e_{kl} > 0$ for an arbitrary non-zero rank-2 tensor $e_{kl}$, respectively.

In this notation equation~\eqref{anisorel1} assumes the form
\begin{equation}\label{anisorel}
\N \, \tilde{\C}(p) + \tilde{\F}(p)\, \tilde{\C}(p) = p^{-2} \, \I
\end{equation}

\begin{theorem} \label{xixi}
1. If the matrix $\N$ is symmetric positive semi-definite and $\F$ is an LICM matrix-valued function satisfying the following conditions: \\ \vspace{0.2cm}
\parbox{\textwidth}{$(\ast\ast_1)$ for every non-zero vector $\V \in \mathbb{R}^N$  the function $\V^\mathsf{T}\, \F(t)\, \V$ does not
vanish identically,}\\ \vspace{0.2cm}
\parbox{\textwidth}{$(\ast\ast_2)$ $\mathbf{F}_0 := \lim_{t\rightarrow 0}\, \mathbf{F}(t)^{-1}$ exists
and the symmetry relations \eqref{symmetries}$_1$ are satisfied,}
then equation~\eqref{convVE} has a single solution $C_{ijkl}(t)$. The function $C_{ijkl}(t)$
is a rank-4 tensor valued Bernstein function satisfying \eqref{symmetries}.\\
2. If $C_{ijkl}(t)$ is a rank-4 tensor-valued Bernstein function such that \\ \vspace{0.2cm}
\parbox{\textwidth}{$(\ast\ast_3)$ for every non-zero symmetric rank-2 tensor $e_{kl}$  the function $C_{ijkl}(t) \, e_{ij}\, e_{kl}$ 
is not constant} and $C_{ijkl}$ satisfies the symmetry relations \eqref{symmetries}$_2$, 
then equation~\eqref{convVE} has a unique solution $\left(\N, \F\right)$. The function $\F$ 
is a rank-4 tensor valued LICM function sa\-tisfying \eqref{symmetries}, while $\N$ is a rank-4 tensor satisfying \eqref{symmetries} 
and $\N \geq 0$ for every rank-2 tensor $e_{kl}$.

\end{theorem}

The first part of the theorem follows from Corollary~\ref{cor4}.\\
The second part follows from Corollary~\ref{cor5}.

\begin{theorem}
Under the hypotheses of Theorem~\ref{xixi} the following relations hold
\begin{enumerate}
\item If $\mathbf{N} > 0$, then $\mathbf{C}(0) = 0$ and $\lim_{t\rightarrow 0} \mathbf{C}^\prime(t) = \mathbf{N}^{-1}$;
\item If $\mathbf{N} = 0$, then $\mathbf{C}(0) = \F_0$.
\item The limit $\lim_{t\rightarrow \infty} \F(t) =: \F_\infty$ always exists and is positive-semi-definite. 
\item If $\F_\infty$ is invertible, then $\lim_{t\rightarrow \infty} \mathbf{C}(t)$ exists and 
\begin{equation} \label{eq3}
\lim_{t\rightarrow \infty} \mathbf{C}(t) = \F_\infty^{\;-1}
\end{equation}
\end{enumerate}

\end{theorem}
\noindent\textbf{Proof.}\\
$$\tilde{\R}(p) = \N + \int_{[0,\infty[} (r + p)^{-1} \, \G(r)\, \mu(\dd r),$$ 
where the Borel measure $\mu$ satisfies equation~\eqref{ineq} and $\vert \G(r) \vert \leq 1$ $\mu$-almost everywhere,
and $\lim_{t\rightarrow 0} \C(t) = \lim_{p\rightarrow \infty} \left[ p^{-1}\, \tilde{\R}(p)^{-1}\right]$, hence 
\begin{equation} \label{Ct}
\lim_{t\rightarrow 0} \mathbf{C}(t) =  \lim_{p\rightarrow \infty} \left\{ p^{-1} \, \left[ \mathbf{N} + 
\int_{[0,\infty[} (r + p)^{-1} \, \mathbf{G}(r)\, \mu(\dd r) \right]^{-1} \right\}
\end{equation}
The second term in the square brackets is the Laplace transform of an LICM function $\F$, defined by equation~\eqref{FLICM}.
 By Proposition~\ref{decrease}  $$\lim_{p\rightarrow \infty} \left[\mathbf{N} +
\int_{[0,\infty[} (r + p)^{-1} \, \mathbf{G}(r)\, \mu(\dd r)\right]  = \mathbf{N}.$$
Thus if $\mathbf{N}$ is invertible, then
$\lim_{p\rightarrow\infty} \left[\mathbf{N} +
\int_{[0,\infty[} (r + p)^{-1} \, \mathbf{G}(r)\, \mu(\dd r)\right]^{-1}$ exists and equals $\mathbf{N}^{-1}$.
It follows from equation~\eqref{Ct} that $\lim_{t\rightarrow 0} \mathbf{C}(t) = 0$ in this case. 

Furthermore, under the same assumption 
$$\lim_{t\rightarrow 0} \mathbf{C}^\prime(t)  = \lim_{p\rightarrow \infty} [p^2\, \tilde{\mathbf{C}}(p)] =
\lim_{p\rightarrow \infty} \left[\mathbf{N} + 
\int_{]0,\infty[} (r + p)^{-1} \, \mathbf{G}(r)\, \mu(\dd r)\right]^{-1} = \mathbf{N}^{-1}$$

If $\N = 0$ then
$$\lim_{t\rightarrow 0} \mathbf{C}(t) = \lim_{p\rightarrow \infty} [p\, \tilde{\mathbf{C}}(p)] = 
\lim_{p\rightarrow \infty} \left[ p \, \tilde{\F}(p)\right]^{-1}$$ 
The limit of the expression on the right-hand side is $\lim_{t\rightarrow 0} \F(t)^{-1} = \F_0$.

Finally, 
$$\lim_{t\rightarrow \infty} \mathbf{C}(t) = \lim_{p\rightarrow 0} [p \, \tilde{\mathbf{C}}(p)] = 
\lim_{p \rightarrow 0} \left[ p \, \mathbf{N} + p \int_{[0,\infty[} (p + r)^{-1} \mathbf{G}(r)\, \mu(\dd r)\right]^{-1}$$
If the limit $\lim_{p \rightarrow 0} \left[ p \int_{[0,\infty[} (p + r)^{-1} \mathbf{G}(r)\, \mu(\dd r)\right]
= \lim_{t\rightarrow \infty} \F(t)$ exists and is invertible 
then equation~\eqref{eq3} is satisfied.

\mbox{ }\hfill $\Box$

We now examine the reverse direction from creep tests to some parameters of the relaxation function. 

\begin{proposition}
1. If $\C(t) \not\equiv 0$ then the function $\C(t)^{-1}$ has a limit  $\G_1 \geq 0$ for $t \rightarrow \infty$.\\
2. If $\G_1 > 0$, then $\lim_{t \rightarrow \infty} \C(t) = \G_1^{\;-1}$.\\
3. The derivative $\C^\prime(t)$ has a limit $\G_2$ for $t \rightarrow \infty$. The limit 
$\lim_{t\rightarrow 0} \C^\prime(t)$ exists, possibly infinite.\\
4. If a $\C(t)$ has a finite limit for $t \rightarrow \infty$, then $\G_2 = 0$.
\end{proposition}
\noindent\textbf{Proof.}\\
Ad 1. The function $\C(t)^{-1}$ is positive semi-definite and non-increasing, hence the limit $\G_1$ exits and is 
positive semi-definite.\\
Ad 2. It follows from the continuity of algebraic inverse at invertible matrices.\\
Ad 3. The function $\C^\prime$ is LICM, hence it is positive semi-definite and non-increasing.  \\
Ad 4. Suppose that $\G_2 > 0$ and $\C(t)$ has a finite limit $\C_\infty$ for $t \rightarrow \infty$. $\C^\prime$ is
non-increasing, hence $\C(t) \geq \C(1) + (t - 1)\, \G_2$. For sufficiently large $t$ the value of $\C(t)$ is larger than $\C_\infty$,
which leads to a contradiction.\\

\mbox{}\hfill $\Box$

\begin{theorem}
If $\mathbf{C}$ is a $\mathcal{S}_+$-valued Bernstein function, 
and \\ \vspace{0.2cm}
\parbox{\textwidth}{$(\ast\ast\ast)$\; for each non-zero vector $\mathbf{v} \in \mathbb{R}^6$ the function $C_{IJ}(t)\, v_I\, v_J$ is not identically zero,} \\
then there is an $\mathcal{S}_+$-valued LICM function $\mathbf{F}$ 
and $\mathbf{N} \in \mathcal{S}_+$ such that the pair $(\N, \F))$
satisfies equation~\eqref{anisorel-t}. 

 Moreover\\
1. If $\C(0) > 0$ then $\N = 0$ and $\lim_{t\rightarrow 0} \F(0) = \C(0)^{-1}$.\\
2. If $\C(0) = 0$, while $\C^\prime_0 := \lim_{t \rightarrow 0} \C^\prime(t)$ is finite and has an inverse then
$\N = \C^\prime(0)^{-1}$.\\
3. $$ \lim_{t\rightarrow \infty} \F(t) = \left\{\begin{array}{cl} 0 & \mbox{if $\lim_{t\rightarrow \infty} \C^\prime(t) > 0$} \\
 \lim_{t \rightarrow \infty} \C(t)^{-1} & \mbox{if $\lim_{t\rightarrow \infty} \C^\prime(t) = 0$}
\end{array} \right.$$
\end{theorem}
\noindent\textbf{Proof.}\\
The derivative $\C^\prime(t)$ is an LICM function and $\C^\prime(t) = \B + \mathbf{Q}(t)$, where $\B = \lim_{t\rightarrow \infty}
\C^\prime(t)$, $\mathbf{Q}$ is LICM  and $\lim_{t \rightarrow \infty} \mathbf{Q}(t) = 0$. Hence
$$\C(t) = \A + t \, \B + \int_0^t \mathbf{Q}(s)\, \dd s. $$
where the matrices $\A$ and $\B$ are symmetric and positive semi-definite. We note that $\A = \C(0)$ and, by Proposition~\ref{decrease}, 
$\B = \lim_{t\rightarrow\infty} \C^\prime(t)$. 

In view of Assumption ($\ast\ast\ast$) the matrix $\tilde{\C}(p)$ has an inverse for every $p > 0$. Since 
$$p \, \tilde{\C}(p) = \A + p^{-1}\, \B + \tilde{\mathbf{Q}}(p) $$ is a Stieltjes function, its algebraic inverse 
$p \, \tilde{\R}(p)$ is a CBF. Hence
$$p \, \tilde{\R}(p) = p \, \N + p \int_{[0,\infty[} (p + r)^{-1} \, \H(r) \, \mu(\dd r)$$
where $\N$ is symmetric positive semi-definite $\mu$ is a Borel measure satisfying \eqref{ineq} and $\H(r)$ a bounded symmetric 
function for $\mu$-almost all $r \in [0,\infty[$. 

By Proposition~\ref{decrease} $\N = \lim_{p \rightarrow \infty} \tilde{\R}(p)$. On the other hand
\begin{multline}
\lim_{p \rightarrow \infty} \tilde{\R}(p) = \lim_{p \rightarrow \infty} \left[p^2\, \tilde{\C}(p)\right]^{-1} 
= \lim_{p \rightarrow \infty} \left[ p \, \A + \B + p\, \tilde{\mathbf{Q}}(p) \right]^{-1} = \\ =
\left\{ \begin{array}{cl} 0 & \A > 0 \\ \left[ \B + \mathbf{Q}(0)\right]^{-1} & \A = 0, \left[ \B + \mathbf{Q}(0)\right] > 0\end{array} \right.
\end{multline}
where we note that $\B + \mathbf{Q}(0) = \C^\prime(0)$. 

Finally
$$
\lim_{t\rightarrow\infty} \F(t) = \lim_{p\rightarrow \infty} \left[ p\, \tilde{\R}(p) \right] =
\left\{ \begin{array}{cl} 0 & \B > 0\\
\left[ \A + \int_0^\infty \mathbf{Q}(t) \, \dd t\right]^{-1} & \B = 0 \end{array}\right.
$$
where we note that $\A + \int_0^\infty \mathbf{Q}(t) \, \dd t = \lim_{t\rightarrow \infty} \C(t)$ if the limit on the right-hand 
side exists.\\
\mbox{ } \hfill$\Box$

Theorem~\ref{xixi} implies that setting the Newtonian viscosity coefficient $\N = 0$ has profound consequences for the creep 
$\C(t)$: the creep either starts with a jump or vertically from the zero value. On the other hand the value of $\N$ can be 
estimated from creep tests as the inverse of the original creep rate. 

In comparison with \cite{HanDuality} the results of Section~\ref{linear} have demonstrated the inseparability of the Newtonian 
component of viscoelasticity from viscoelastic memory effects. 

\section{Conclusions.}

We have demonstrated a particular role of LICM kernels in two classes of convolution equations and utility
of the concepts of CBFs and Stieltjes derivatives in the study of existence problems for these equations. A particular class 
of convolution equations studied here are fundamental in linear viscoelasticity. 

Another convolution equation 
has allowed us to define a class of anisotropic generalized fractional derivatives associated with matrix-valued LICM kernels. We have also shown 
that the kernel appearing in the definition of anisotropic GFD and the associated anisotropic fractional integrals should be singular at 0.  

\appendix

\section{A remark on the convolution algebra.}

For our purposes it is important that the convolution algebra has a unity. The unity is not the convolution with a function.
Hence the convolution algebra must include Borel measures. The convolution $\rho\ast\nu$ of two measures $\rho$ and $\nu$ 
defined on $[0,\infty[$ is
defined as the Borel measure $\lambda$ satisfying the identity 
$$\int_{[0,\infty[} f(r)\, \dd \lambda(\dd r) = \int_{[0,\infty[} \int_{[0,\infty[} f(r + s)\,\rho(\dd r)\,\nu(\dd s)$$
for every continuous function $f$ with compact support.

This definition is easily extended to matrix-valued measures.

For our purposes the convolution algebra has to involve only Borel measures of the form $u \, \C + \F(t) \, \dd t$,
where the unity is a measure defined in Section~\ref{First}.

\section{Matrix-valued Stieltjes functions and Complete Bernstein functions.}\label{appAniso}

We shall now use some results from Appendix~B of \cite{HanAnisoWaves}.

A matrix-valued Stieltjes function $\mathbf{Y}(p)$ has the following integral representation:
\begin{eqnarray}\label{mStieltjes}
\mathbf{Y}(p) = \mathbf{B} + \int_{[0,\infty[} (p + r)^{-1} \, \mathbf{H}(r) \, \mu(\dd r) 
\end{eqnarray}
where $\mathbf{B} \in \mathcal{S}_+$, $\mu$ is a Borel measure on $[0,\infty[$ satisfying \eqref{ineq} and 
 $\mathbf{H}(r)$ is an $\mathcal{S}_+$-valued function defined and bounded 
$\mu$-almost everywhere on $[0,\infty[$. \\
Conversely, any matrix-valued function with the integral representation 
\eqref{mStieltjes} is an $\mathcal{S}_+$-valued Stieltjes function.

\begin{theorem}
Every matrix-valued Stieltjes function is the Laplace transform of an element of $\mathcal{E}$.
\end{theorem}

\noindent\textbf{Proof.}\\

The Laplace transform of the $\mathcal{S}_+$-valued LICM function $\mathbf{A}(t)$ (equation~\eqref{LICM}) 
is given by the equation 
\begin{equation} \label{mlapLICM}
\tilde{\mathbf{F}}(p) = \int_{[0,\infty[} (p + r)^{-1} \, \mathbf{H}(r) \, \mu(\dd r)
\end{equation}
where $\mu$, $\H$ satisfy the same conditions as in \eqref{mStieltjes}.

The second term on the right-hand side of equation~\eqref{mStieltjes} involves a double Laplace transform,
hence it equals 
$$\mathbf{V}(t) := \int_0^\infty \e^{-p t} \left[\int_0^\infty \e^{-r t} \, \mathbf{H}(r) \, \mu(\dd r) \right] \dd t$$
where the Borel measure $\mu$ satisfies inequality~\eqref{ineq}. The inner integral represents a general matrix-valued LICM $\F(t)$.
Thus $\mathbf{V}(t)$ is the Laplace transform of a general matrix-valued LICM and $\mathbf{Y}(p)$ is the Laplace transform of a general 
element of $\mathcal{E}$.\\

\mbox{ } \hfill$\Box$

An $\mathcal{S}_+$-valued CBF $\mathbf{Z}(p)$ has the following integral representation:
\begin{eqnarray}\label{mCBF}
\mathbf{Z}(p) = p\, \mathbf{B} + p \int_{[0,\infty[} (p + r)^{-1} \, \mathbf{H}(r) \, \nu(\dd r) 
\end{eqnarray}
where $\mathbf{B} \in \mathcal{S}_+$, $\nu$ is a Borel measure on $[0,\infty[$ satisfying \eqref{ineq} and 
$\mathbf{H}(r)$ is an $\mathcal{S}_+$-valued function defined
$\nu$-almost everywhere on $[0,\infty[$.

Conversely, any $\mathcal{S}_+$-valued function with the integral representation 
\eqref{mCBF} is a $\mathcal{S}_+$-valued CBF.

It follows immediately that the the function $p^{-1} \, \mathbf{Z}(p)$, where $\mathbf{Z}$ is an $\mathcal{S}_+$-valued 
CBF function, is an $\mathcal{S}_+$-valued Stieltjes function. 

We quote Lemma~1 in Appendix~B of op. cit. in the form of the following theorem 
\begin{theorem} \label{inv}
If $\mathbf{Z}(p)$ is an $\mathcal{S}_+$-valued CBF and does not vanish identically, then $\mathbf{Z}(p)^{-1}$ is an 
$\mathcal{S}_+$-valued Stieltjes function.

Conversely, if $\mathbf{Y}(p)$ is an $\mathcal{S}_+$-valued function does not
vanish identically then $\mathbf{Y}(p)^{-1}$ is a CBF.

\end{theorem}

\providecommand{\bysame}{\leavevmode\hbox to3em{\hrulefill}\thinspace}
\providecommand{\MR}{\relax\ifhmode\unskip\space\fi MR }
% \MRhref is called by the amsart/book/proc definition of \MR.
\providecommand{\MRhref}[2]{%
  \href{http://www.ams.org/mathscinet-getitem?mr=#1}{#2}
}
\providecommand{\href}[2]{#2}

\end{document}